\documentclass[12pt,letterpaper,reqno]{amsart}

\usepackage{times} 
\usepackage[T1]{fontenc} 
\usepackage{mathrsfs} 
\usepackage{latexsym} 
\usepackage[dvips]{graphics}
\usepackage{epsfig}
\usepackage{amsmath,amsfonts,amsthm,amssymb,amscd} 
\input amssym.def 
\input amssym.tex

\newcommand\be{\begin{equation}} 
\newcommand\ee{\end{equation}}
\newcommand\bea{\begin{eqnarray}} 
\newcommand\eea{\end{eqnarray}} 
\newcommand\bi{\begin{itemize}}
\newcommand\ei{\end{itemize}} 
\newcommand\ben{\begin{enumerate}} 
\newcommand\een{\end{enumerate}}
\newcommand\bc{\begin{center}} 
\newcommand\ec{\end{center}} 
\newcommand\ba{\begin{array}} 
\newcommand\ea{\end{array}}

\theoremstyle{definition} 

\begin{document}

\title[Why should one expect to find long runs of 
(non)-Ramanujan primes ?] 
{Why should one expect to find long runs of (non)-Ramanujan primes ?}


\author{Peter Hegarty} \address{Department of Mathematical Sciences, 
Chalmers University Of Technology and University of Gothenburg,
41296 Gothenburg, Sweden} \email{hegarty@chalmers.se}


 \keywords{}

\date{\today}

\begin{abstract} 
Sondow et al have studied Ramanujan primes (RPs) and observed numerically that, while half of all primes are RPs asymptotically, one obtains runs of consecutives RPs (resp. non-RPs) which are statistically significantly longer than one would expect if one was tossing an unbiased coin. In this discussion paper we attempt a heuristic explanation of this phenomenon. Our heuristic follows naturally from the Prime Number Theorem, but seems to be only partly satisfactory. It motivates why one should obtain long runs of both RPs and non-RPs, and also longer runs of non-RPs than of RPs. However, it also suggests that one should obtain longer runs of RPs than have so far been observed in the data, and this issue remains puzzling.
\end{abstract}


\maketitle

\setcounter{equation}{0}

\setcounter{equation}{0}

\section{The model}

Consider the following random process : you have an infinite supply of 
identical biased 
coins, which 
return heads with probability $p \in (\frac{1}{2},1]$, and tails with 
probability $1-p$. Now toss $N$ of these coins. For each $i = 1,...,N$, let 
$h_i, t_i$ denote the number of heads (resp. tails) among the first
$i$ tosses. Thus $h_i + t_i = i$. Also denote 
\be
\Delta_i := h_i - t_i.
\ee
The coins that come up as heads will each be colored red or blue, according
to the following rule : Suppose the $i$:th toss is a head. Then color this
coin red if and only if the following two conditions sare satisfied : 
\be
\Delta_j \geq \Delta_i, \;\; {\hbox{for all $j = i+1,...,N$}}
\ee
and
\be
{\hbox{If $1 \leq k < i$ and $\Delta_k \geq \Delta_i$, then there exists
$l \in (k,i)$ such that $\Delta_l < \Delta_i$.}}
\ee
Otherwise, color the coin blue. The definition requires some thought, so here is
an example to illustrate how the scheme works :

\begin{table}[ht!]
\begin{center}
\begin{tabular}{|c|c|c|c|c|c|c|c|c|c|c|c|c|c|c|c|c|} \hline
${\hbox{toss}}$ & $1$ & $2$ & $3$ & $4$ & $5$ & $6$ & $7$ & $8$ & $9$ & $10$ & $11$ & $12$ & $13$ & $14$ & $15$ & $16$ \\ \hline
${\hbox{result}}$ & $H$ & $T$ & $H$ & $H$ & $T$ & $H$ & $H$ & $H$ & $T$ & $H$ & $T$ & $H$ & $H$ & $T$ & $H$ & $H$ \\ \hline
$\Delta$ & $1$ & $0$ & $1$ & $2$ & $1$ & $2$ & $3$ & $4$ & $3$ & $4$ & $3$ & $4$ & $5$ & $4$ & $5$ & $6$ \\ \hline
${\hbox{color}}$ & $B$ & $-$ & $R$ & $B$ & $-$ & $R$ & $R$ & $B$ & $-$ & $B$ & $-$ & $R$ & $B$ & $-$ & $R$ & $R$ \\ \hline
\end{tabular}
\end{center}
\end{table}    
$\;$ \par
\par
Now what do we expect to observe, when $N$ is large ? Well, with high
probability (w.h.p.), we will observe close to $pN$ heads and close to 
$(1-p)N$ tails. Hence $\Delta_N \approx (2p-1)N$ w.h.p. Now at most one coin is
colored red for each positive value attained by the function $\Delta$. On the
other hand, for any $\epsilon > 0$, w.h.p. the function $\Delta$ will
eventually exceed $(2p-1-\epsilon)N$ for good. Hence, w.h.p. about $(2p-1)N$ 
coins will be colored red, and the remaining heads, about $(1-p)N$ in number,
will be colored blue. So, as $N \rightarrow \infty$, the fraction of
redheads, amongst all heads, will almost surely (a.s.) approach 
$\frac{2p-1}{p}$ and the fraction
of blueheads will a.s. approach $\frac{1-p}{p}$. In particular, when 
$p = 2/3$, about half the head-coins will be colored red and half colored blue. 
\\
\\
I claim that, with $p = 2/3$, this is a good basic model to have in mind when 
one considers
Ramanujan primes (RPs) : the red coins corresponding to RPs and the 
blue ones to non-RPs. I will explain two things :
\\
\\
1. Why this model is reasonable.
\\
2. Why one expects to get longer runs of blue coins than if
the red-blue coloring was done by tossing another, fair coin. 
\\
\\
I will deal with the second issue first. However, the analysis will show that,
in this model, one also expects longer runs of red coins than
if the coloring was done fairly at random, though not as long as the 
blue runs. This may seem to contradict
the data in \cite{SNN}. 
After explaining why I nevertheless consider the model to be
reasonable, I will discuss this issue. 
     
\setcounter{equation}{0}

\section{Why do we get long monochromatic runs ?}

If a biased coin with probability $p$ of heads is tossed $N$ times, then it is
well-known that the expected length of the longest run of consecutive
heads is approximately $\frac{\log N}{\log (1/p)} = \log_{1/p} N$. Another
way of looking at this is that, for any $k \in \mathbb{N}$, one expects to 
have to toss 
the coin on the order of $(1/p)^k$ times to have a reasonable probability of 
seeing at least one run of $k$ consecutive heads. This is easy to see
intuitively : the probability of any $k$ consecutive coin tosses all
resulting in heads is $p^k$ and thus, by linearity of expectation, the 
expected number of such runs amongst $N$ tosses is $(N-k+1)p^k$, which 
(for any fixed $k$) will
be $\Theta(1)$ when $N = \Theta[(1/p)^k]$. 
\par In particular, when $p = 1/2$, we expect to have to make on the order
of $2^k$ tosses to have a reasonable probability of witnessing a run of 
$k$ heads. 
\\
\\
The following facts about biased coin-tossing are also well-known :
\\
\\
{\bf Proposition 2.1} {\em Suppose we toss a sequence of identical biased
coins with probability $p > 1/2$ of heads. With notation as in Section 1, 
for each $N \in \mathbb{N}$, let
$c_{N,p}$ denote the probability that $\Delta_i \geq 0$ for all
$i = 1,...,N$. Then the numbers $c_{N,p}$ are non-increasing in $N$ and if
we let $c_p := \lim_{N \rightarrow \infty} c_{N,p}$, then  
\be
c_{p} = \frac{2p-1}{p} > 0.
\ee}
{\sc Proof} : That $c_{N,p} \geq c_{N+1,p}$ is trivial. Let $E,F$ and 
$F^{\prime}$ denote
the following three events :
\par $E$ : the event that $\Delta_i > 0$ for all $i > 0$, when we make
an infinite sequence of tosses.
\par $F$ : the event that $\Delta_i \geq 0$ for all $i > 0$ when we make an 
infinite sequence of tosses.
\par $F^{\prime}$ : the event that $\Delta_i \geq \Delta_1$ for all $i > 1$
when we make an infinite sequence of tosses.
\\
\\
Since the coin tosses are independent, one has 
\be
\mathbb{P}(F) = \mathbb{P}(F^{\prime}). 
\ee
Secondly, it is clear that
\be 
c_{p} = \mathbb{P}(F).
\ee 
Thirdly, the
event $E$ occurs if and only if the first toss yields a head and thereafter
event $F^{\prime}$ occurs. Hence,
\be
\mathbb{P}(E) = p \cdot \mathbb{P}(F^{\prime}).
\ee
In Example 1.13, Chapter 3 of \cite{D}, it is shown that 
\be
\mathbb{P}(E) = 2p-1.
\ee
Eqs. (2.2)-(2.5) together imply (2.1), and the proof is complete.
\\
\\
{\bf Remark 2.2} It is no accident that $c_p$ equals the fraction of 
redheads in the head-coloring model of Section 1. Indeed, this observation
is the basis for the rigorous proof of (2.1). 
\\
\\
Fix $p > 1/2$ and consider the head-coloring 
model of Section 1. Fix $k \in \mathbb{N}$ and let
$\mathcal{E}_{k,N}$ denote the expected number of runs of $k$ 
redheads, when $N$ coins are tossed. I claim that, as $N \rightarrow \infty$, 
\be
\left[ \frac{(2p-1)^2}{p} \right] \cdot p^{k-1} \lesssim 
\frac{\mathcal{E}_{k,N}}{N} \leq 
p^{k-1}.
\ee
To see the right-hand inequality, just observe that if $k$ consecutive heads
are all colored red, then at the very least the $k-1$ heads from the 2nd to the
last must have 
been a run of $k-1$ heads, with no tails in between. This
happens with probability $p^{k-1}$. Thus, the probability of a run 
of $k$ redheads with a fixed starting point is at most  
$p^{k-1}$. Since there are $N-k+1$ possible starting points
for the run, linearity of expectation implies that 
\be
\mathcal{E}_{k,N} \leq (N-k+1) p^{k-1},
\ee
which gives the right-hand inequality in (2.6). For the lower bound, we
again consider a fixed starting point. A sufficient condition to get
a run of $k$ redheads with a given starting point is that the following
three events all occur :
\par $A$ : the starting point is a redhead,
\par $B$ : it is followed by a run of $k-1$ heads, with no tails in between
\par $C$ : the value of $\Delta$ never again goes below its value at the 
end of this run of $k$ heads.
\\
\\
It is clear that each of $B$ and $C$ positively correlates with $A$, while $B$ 
and $C$ are independent of one another. Hence 
\be
\mathbb{P}(A \wedge B \wedge C) \geq \mathbb{P}(A) \times \mathbb{P}(B) \times
\mathbb{P}(C).
\ee
As shown in Section 1, we know that $\mathbb{P}(A) \rightarrow 2p-1$ as
$N \rightarrow \infty$. As above, the event $B$ occurs with 
probability $p^{k-1}$. Thirdly, it is immediate that $\mathbb{P}(C) \geq c_p
= \frac{2p-1}{p}$. Plugging everything into (2.8), we find that
the probability of a run of $k$ redheads with a given starting point is
at least $\left[ \frac{(2p-1)^2}{p} \right] \cdot p^{k-1}$. Linearity
of expectation then yields the left-hand inequality in (2.6). 
\par This brings us to our first result :
\\
\\
{\bf Proposition 2.3} {\em In the model of Section 1, the expected length
of the longest run of consecutive reds among a total of $N$ heads, is
on the order of $\frac{\log N}{\log (1/p)}$. In other words, for any 
$k \in \mathbb{N}$, we expect to have to make on the order of
$p^k$ tosses in order to have a reasonable
probability of observing a run of $k$ redheads.
\par In particular, when $p = 2/3$, the expected length
of the longest run of consecutive reds among a total of $N$ heads, is
on the order of $\frac{\log N}{\log (3/2)}$. In other words, for any 
$k \in \mathbb{N}$, we expect to have to make on the order of
$\left( \frac{3}{2} \right)^k$ tosses in order to have a reasonable
probability of observing a run of $k$ redheads.}
\\
\\
{\bf Remark 2.4} The proposition says that, for any fixed $p > 1/2$ and 
very large $N$, we expect to see runs of redheads amongst the heads 
of similar length to runs of heads amongst all the coins.
\\
\\ 
So what about blues ? Here, for simplicity, I only consider the case $p = 2/3$
for the moment{\footnote{I will generalise to arbitrary $p > 1/2$ when I get 
the time. I will indicate below where changes need to be made.}}. 
Let $k \in \mathbb{N}$ and let 
$\mathcal{F}_{2k,N}$ denote the 
expected number of runs of $2k$ consecutive blue heads somewhere 
amongst the first $N$ coins. Suppose, for example, that in 
a run of $3k$ consecutive tosses one observes at least $2k$ 
tails{\footnote{More generally, one will need to replace 2 and 3 by some numbers depending on $p$, and chosen in such a way that the final exponent
in (2.13) will be less than $p$.}}. This means 
that the function $\Delta$ will have decreased by at least $k$ over this run.
All succeeding heads will definitely 
be colored blue at least until the function $\Delta$
has risen by $k$ again. 
Then it is a tedious, but standard, calculation to show that there is a
fixed $u > 0${\footnote{This number will also 
depend on $p$ in a general analysis.}} such that, in order for the function 
$\Delta$ to increase by $k$, at least $2k$ heads will need to 
be revealed.  
Hence, as 
$N \rightarrow \infty$, 
\be
\frac{\mathcal{F}_{2k,N}}{N} \gtrsim u \cdot q_{2k},
\ee
where $q_{2k}$ is the probability of a run of $3k$ tosses yielding at least 
$2k$ tails. Explicitly, one has 
\be
q_{2k} = \sum_{l=2k}^{3k} \left( \begin{array}{c} 3k \\ l \end{array} \right)
\left( \frac{1}{3} \right)^l \left( \frac{2}{3} \right)^{3k-l}.
\ee
A lower bound for this is got by simply taking the $l = 2k$ term, hence 
\be
q_{2k} \geq \left( \begin{array}{c} 3k \\ 2k \end{array} \right)
\left( \frac{1}{3} \right)^{2k} \left( \frac{2}{3} \right)^{k}.  
\ee
Now let $k \rightarrow \infty$ also. Put $v_k := \sqrt{\frac{3}{4\pi k}}$.
Applying Stirling's estimate to (2.11),
one easily computes that
\be
\left( \begin{array}{c} 3k \\ 2k \end{array} \right) \sim 
v_k \left( 
\frac{3^3}{2^2} \right)^k = v_k \left( \frac{27}{4} \right)^k,
\ee
and hence that  
\be
q_{2k} \gtrsim v_k \left( \frac{1}{2} \right)^k = v_k \left( \frac{1}{\sqrt{2}} 
\right)^{2k}.
\ee
Putting all this together, and using the fact that $u$ in (2.9) is a 
constant, plus that the function $v_k$ 
decreases subexponentially in $k$, we have our second result :
\\
\\
{\bf Proposition 2.5} {\em In the model of Section 1, with 
$p = 2/3$, the expected length
of the longest run of consecutive blues among a total of $N$ heads, is
at least $\frac{\log N}{\log (\sqrt{2})} = 2 \log_{2} N$. In other words, 
for large but fixed
$k \in \mathbb{N}$, we expect to have to make at most
on the order of $(\sqrt{2})^k = 2^{k/2}$ tosses in order to have a reasonable
probability of observing a run of $k$ blueheads.} 
\\
\\
{\bf Remark} Note
the use of the words $\lq$at least' and $\lq$at most' in Proposition 2.5, 
as against $\lq$approximately' in Proposition 2.3. 
This reflects the fact that we only have a lower bound in (2.9), whereas
in (2.6) we have both upper and lower bounds. While it seems
difficult to compute $\frac{\mathcal{F}_{2k,N}}{N}$ exactly, it is 
quite easy to see, with the help of Proposition 2.1, 
that there will be some upper bound of the form 
$c_{1}^{2k}$, for some $c_1 < 1$.  
Hence, in order to have a 
reasonable probability of observing a run of $k$ blueheads, one does 
expect to have to make a number of tosses which is exponential in $k$.  
\\
\\
From Propositions 2.3 and 2.5 it follows that one expects to see longer
runs, both of reds and blues, than if the coloring was done fairly at
random, but that one expects to see even longer runs of blueheads than of
redheads. Or, to put it another way, for any fixed, and large enough 
$k \in \mathbb{N}$, 
one expects to see a run of $k$ blueheads somewhat earlier
than a run of $k$ redheads, and one expects to see both in turn much 
earlier than if the coloring was done fairly at random.  
       
\setcounter{equation}{0}

\section{Why is this a reasonable model for RPs ?}

First of all, we just focus on the $\lq$ordinary' RPs discussed in 
\cite{SNN}{\footnote{Table 1 of this paper contains some errors, but these have been corrected in arXiv:1105.2249(v2)}}, which we wish to compare with the $p = 2/3$ model above. 
The obvious 
extension to the so-called Generalised RPs of \cite{ABMRS} will be 
mentioned at the end. 
\\
\\
So consider $p = 2/3$ as fixed for now. Let me describe, in somewhat informal 
terms, 
another random model which I claim is
asymptotically equivalent to that in Section 1. Suppose we have
two different radioactive substances, H (head) and T (tail). H decays
twice as fast as T (i.e.: T has double the half-life of H). At 
some time $t = 0$, I start observing both substances and record each
individual decay of a H- or a T-atom. For each $i \in \mathbb{N}$, let 
$h_i, t_i$ denote the number of H- (resp. T-) atoms amongst the 
first $i$ which decay. Then color a decayed H-atom red if conditions
(1.2) and (1.3) hold, otherwise blue.
\par If one wants to be more formal, one can phrase this in terms of
two independent, parallel Poisson processes, one of which has double
the intensity of the other. But the point is that this model is
equivalent, as $t \rightarrow \infty$, with that of Section 1. I leave it 
to the reader to convince himself of this. 
\\
\\
From here, we can see the relevance to Ramanujan primes. First of all, the
prime number theorem says that, for large $x$, there are approximately
$\frac{x}{\log x}$ primes up to $x$. Equivalently, it says that, for large
$n$, the $n$:th prime satsifies $p_n \sim n \log n$. Now suppose that
we start from some $p_n$ and begin searching to the right of it for the 
next prime $p_{n+1}$. There is, of course, nothing random about this 
search process. But in a well-known $\lq$random model' for prime gaps, one
imagines that this search is a Poisson process with intensity 
$\frac{1}{\log n} \sim \frac{1}{\log p_n}$. In this random model, there
is also no need to start the search at a prime : the starting point can also 
be chosen at random without affecting the model.  
\par Now suppose one chooses a large random number $x$ and 
starts two prime searches in parallel, one at $x$ and the other at $x/2$, where
the former search proceeds twice as fast as the latter. By this I mean that,
when in the former process we have searched from $x$ up to $x+t$, then
in the latter we have searched from $x/2$ up to $(x+t)/2$.  
Since $\log (x/2) \sim \log x$, one sees that the former Poisson process
has approximately twice the intensity of the latter. Clearly, the method of
determining whether a prime is Ramanujan or not is basically the same as 
that of deciding how to color the primes $\lq$revealed' in the former of these
two Poisson processes. Hence, we have shown how our model in Section 1
is a reasonable model for the Ramanujan primes.   
\\
\\
Finally, we turn to Generalised RPs. To model $c$-Ramanujan primes, one 
should choose $p = \frac{1}{1+c}${\footnote{I will leave all further 
calculations for general $c$ to another forum.}}. 

\setcounter{equation}{0}

\section{Discussion}

Everything above is, of course, heuristics. There is nothing $\lq$random' about 
the prime numbers. More importantly, the Prime Number Theorem is 
a very precise statement about the density of the primes. This suggests that
a main problem with our heuristic model is its Markovian nature. 
In other words, if I start a search for a prime from some point $x$, and
don't find any prime up to $x + t$, say, then the PNT implies that
I am now, in some sense, $\lq$more likely' to find a prime between 
$x+t$ and $x+2t$. The further one searches the 
more restrictive the PNT becomes. Since, in order to observe a 
monochromatic run (equivalently, a run of (non)-RPs) of length $k$, one 
expects to have to search in a range exponential in $k$ (Props. 2.3 and 2.5),
the PNT will, indeed, impose severe restrictions on what can happen. 
\par In \cite{SNN}, the authors fix an upper bound $x$, and find the 
length of the longest run of (non)-RPs in $[1,x]$. Our models suggest that 
a better way of collecting the data would be to fix an integer $k$, and 
determine the first time a run of $k$ (non)-RPs appears. This is 
because, $\lq$locally', the model of two parallel Poisson processes seems
to be okay. Despite the problems discussed above, I do not see clearly why 
Propositions 2.1 and 2.2 should not be a reasonable guide what should be 
observed, as $k \rightarrow \infty$. In other words, I conjecture that,
once the numbers get big enough, one will also observe considerably longer runs
of Ramanujan primes than a fair coin-tossing model would suggest, though
never as long as the runs of non-Ramanujan primes. If this is not 
the case, if the data in \cite{SNN} is a reasonable guide
to what happens asymptotically, then our model must contain a serious flaw. 

\section*{Acknowledgements}

I thank Jonathan Sondow and Johan Tykesson for helpful discussions.

\end{document}